\documentclass[12pt]{article}
\usepackage{amscd, amssymb,latexsym,amsmath, amscd, amsmath}
\usepackage[all]{xy}
\begin{document}

\title{{Reduction of homomorphisms  
mod $p$ and algebraicity} 
\date{}
\author{Chandrashekhar Khare \& Dipendra Prasad}}
\maketitle

\begin{abstract}
{Let $K$ be a number field, and $A_1,A_2$ abelian varieties over $K$.
Let $P$ (resp. $Q$) be a  non-torsion point in $ A_1(K)$ (resp. $A_2(K)$) 
such that for almost
  all places $v$ of $K$, the order of $Q$ mod $v$ divides the order of $P$
  mod $v$. Then we prove (under some conditions on $A_i, i=1,2$) 
that there is 
a homomorphism $j$ from $A_1$ to $A_2$ such that $j(P) = Q$. We formulate and 
extend such a result for any subgroups of
$A_i(K), i=1,2$. These results in particular 
extend the work of Corrales and Schoof
from elliptic curves to abelian varieties.}
\end{abstract}

\newtheorem{theorem}{Theorem}
\newtheorem{lemma}{Lemma}
\newtheorem{prop}{Proposition}
\newtheorem{cor}{Corollary}
\newtheorem{example}{Example}
\newtheorem{conjecture}{Conjecture}
\newtheorem{definition}{Definition}
\newtheorem{question}{Question}
\newtheorem{memo}{Memo}
\newcommand{\rhobar}{\overline{\rho}}
\newcommand{\Sha}{{\rm III}}
\newtheorem{conj}{Conjecture}

\section{Introduction}
 
Let $A_1$ and $A_2$ be two abelian varieties over a number field $K$. 
There is a finite set $S$ of places of  $K$ such that there are
abelian schemes ${\cal A}_i$ over ${\cal O}_S$, 
the ring of $S$-integers of $K$,
with generic fibres $A_i$. For $v$ not in $S$ (all the places $v$ we
consider below are {\it outside}  $S$ even if this is not mentioned
explicitly) we can consider reduction mod $v$ of the schemes ${\cal A}_i$ 
at $v$. We will abuse notation to denote these abelian varieties 
over finite fields also as $A_i$. Consider the specialisation map $sp_v:A_i(K)
\rightarrow A_i(k_v)$, and denote the image of
$A_i(K)$ by $A_i(v)$. 
We say that a homomorphism $\phi:A_1(K) \rightarrow A_2(K)$
specialises mod $v$ if there
is a homomorphism $\phi_v:A_1(v) \rightarrow A_2(v)$ such that the
diagram 
$$
\begin{CD}
  A_1(K) @>\phi>> A_2(K) \\
  @V{sp_v}    VV    @VV{sp_v}   V \\
  A_1(v) @>\phi_v>>  A_2(v)
\end{CD}
$$ commutes. The aim of this paper is to answer the following question.

\begin{question}\label{special}Let $A_1$ and $A_2$ be abelian varieties over 
a number field $K$. 
  Let $\phi:A_1(K) \rightarrow A_2(K)$ be a
  homomorphism of abelian groups that specialises mod $v$ for almost
  all places $v$ of $K$
  to give a homomorphism $\phi_v$ as above.
  Is the restriction of $\phi$ to a subgroup of finite
  index of $A_1(K)$ induced by a homomorphism $\alpha_{\phi} \in   {\rm
  Hom}_K(A_1,A_2)$?
\end{question}

\noindent{\bf Remark:} Although the question above is posed for 
a homomorphism from $A_1(K)$ 
to $A_2(K)$, the heart of the question is ``pointwise'': namely
given points $P$ in $A_1(K)$, and $Q$ in $A_2(K)$ such  that 
for almost all places $v$ of $K$ 
the order of $Q$ in $A_2(k_v)$ divides the order
of $P$ in $A_1(k_v)$, then is $Q$ related to $P$  by a homomorphism from 
$A_1$ to $A_2$? 
(Because of Lemma \ref{abstract} in the next section we see
that this is the essential content of the question.)

Here are the main theorems of this paper.
\begin{theorem}\label{main}
  Let $A_1$ and $A_2$ be two simple 
abelian varieties over a number field $K$ 
with the endomorphism rings over the algebraic closure of $K$ 
as ${\Bbb  Z}$.  Assume that  ${\rm dim}(A_i)=g_i$ is 
either odd or $g_i=2$ or 6. (We do not assume that $\dim A_1 = \dim A_2$.)
Let $P$ (resp. $Q$) be a  non-torsion point in $ A_1(K)$ (resp. $A_2(K)$) 
such that for almost
  all places $v$ of $K$, the order of $Q$ mod $v$ divides the order of $P$
  mod $v$. Then $A_1$ and $A_2$ are isogenous, and there is 
an isogeny $j$ from $A_1$ to $A_2$ such that $j(P) = nQ$ for some 
integer $n$. 
Furthermore, if one of the abelian
varieties is an elliptic curve, we have the same conclusion without
any restriction on its endomorphism ring.

\end{theorem}

\begin{theorem}
  Let $A_1$ and $A_2$ be two simple abelian varieties
 over a number field $K$ with 
${\rm End}_{\bar{K}}(A_i)={\Bbb  Z}$ and such that ${\rm dim}(A_i)=g_i$ is 
either odd or $g_i=2$ or 6. In case one of the $A_i$'s is an elliptic curve,
we allow the elliptic curve to have complex multiplication. Then we have the following.
  \begin{enumerate} \item 
For any abstract  homomorphism $$\phi:A_1(K) \rightarrow A_2(K)$$ that 
specialises mod $v$ for almost all places $v$ of $K$, 
an integral multiple of $\phi$ is given by 
a homomorphism from $A_1$ to $A_2$. 
\item If 
$A_1(K)$ is infinite, $\phi$ itself
  is given by an isogeny. 
\end{enumerate}
\end{theorem}

The restrictions on dimension in the two theorems 
arise from those of Serre's theorems
in article 137 of [S-IV]. We prove  theorem 1 in section 3. 
In section 4 we prove a general lemma that 
answers  Question \ref{special} in the case of CM elliptic curves.

The main component of the proof given in this paper, besides the theorem of Serre
on the image of Galois group on the $\ell$-adic Tate module of an Abelian variety,
is the Kummer theory on Abelian varieties, initiated and studied by Ribet in
several papers, and completed in his paper with Jacquinot, cf. [JR]. 
We recall these theorems from [JR] in section 2. We also
make crucial use of a theorem of Bogomolov, cf. [B], that 
the image of Galois group on the $\ell$-adic Tate module of an 
Abelian variety contains
an open subgroup of the homotheties.

The method we employ in the paper will be able to 
generalise theorem 1 to general abelian varieties
if an analogue of Serre's theorem was available 
for them.  

We note that Theorems 1 and 2 in the special case of
elliptic curves is due to Corrales and Schoof, see [CS]. 
However, our proof for elliptic curves is quite different
from the proof of [CS].

In sections 6 and 7 we look at the analogue of these questions for the
case of linear algebraic groups, and general rational varieties 
respectively.

\section{Theorems due to Jacquinot and Ribet} In this section we recall two theorems
due to Jacquinot and Ribet [JR] which play a fundamental role in this paper.

For any abelian variety $A$ over a number field $K$, we let $A[\ell]$ denote the 
set of $\ell$-torsion  points of $A$ (defined over the algebraic closure $\bar{K}$
of $K$). For any point $P$ in $A(\bar{K})$, $P/\ell$ denotes {\it any} point $Q$ in
$A(\bar{K})$ such that $\ell Q =P$. The minimal field extension 
of $K$ (inside $\bar{K}$) over which
every point of $A[\ell]$ and $Q$ is defined is independent of the choice of
$Q$, and is denoted by $K(A[\ell], P/\ell)$.

The following is theorem 3.1 of [JR].

\begin{theorem} Let $A$ be an abelian variety over a number field $K$, and
let $P \in A(K)$. Let $B$ be the connected component of 0 in the Zariski
closure of the subgroup of $A(K)$ generated by $P$. Then the Galois group of
$K(A[\ell], P/\ell)$ over $K(A[\ell])$ is naturally isomorphic to $B[\ell]$
for almost all primes $\ell$.
\end{theorem}

From this theorem, [JR] deduce the following very important result for us.

\begin{theorem}  Let $A$ be an abelian variety over a number field $K$, and
$P_1,P_2$ be points of $A(K)$. Then there exists a $K$-endomorphism $f:A \rightarrow 
A$, and an integer $n \geq 1$ such that $f(P_1) = nP_2$ if and only if,
$$K(A[\ell], P_2/\ell) \subseteq K(A[\ell], P_1/\ell),$$
for an infinite number of primes $\ell$.
\end{theorem}

\section{Proof of Theorem \ref{main}}

We will give the proof of  theorem 1 when neither of the 
$A_i$'s is an elliptic curve with complex multiplication. The proof when
one of the $A_i$'s is an elliptic curves with complex multiplication is
very similar, and will be left to the reader.

Under our assumptions and by the theorems of Serre in [S-IV], for all
$\ell$ large enough (see Th\'eor\'eme C on page 40 and the Corollaire on
page 51 of the letter to Vign\'eras in [S-IV]), 
the image of the action
of $G_K$ on $A_i[\ell]$ 
is $GSp_{2g_i}({\Bbb Z}/{\ell})$. 

From theorem 3 (due to [JR]) recalled above, we know that 
for almost all primes $\ell$,   
the field extension $K(A_1[\ell],P/{\ell})$ of $K$ obtained by 
attaching the $\ell$-torsion points $A_1[\ell]$ inside $A_1(\bar{K})$,
and any $\ell$-division point $P/{\ell}$ of $P$, is a Galois extension
of $K$ with Galois group $G_{P,\ell}$ which sits in the following 
split exact sequence,
$$0 \rightarrow A[\ell] \rightarrow G_{P,\ell} \rightarrow
GSp_{2g_1}({\Bbb Z}/\ell) \rightarrow 1.$$
We have a similar statement for the  
extension $K(A_2[\ell], Q/{\ell})$ of $K$.

In what follows, we will be considering only those primes $\ell$ for
which Serre's theorem, and theorem 3 (due to  
[JR]) applies for the points $P$ and $Q$.
We will be using below the structure of the Galois groups arising from these 
theorems without explicit mention. In particular, these theorems imply that
$K(A_1[\ell], A_2[\ell], P/\ell)$ is an abelian extension of
$K(A_1[\ell], A_2[\ell])$ with Galois group $({\Bbb Z}/\ell)^{2g_1}$, and  
similarly 
$K(A_1[\ell], A_2[\ell], Q/\ell)$ is an abelian extension of
$K(A_1[\ell], A_2[\ell])$ with Galois group $({\Bbb Z}/\ell)^{2g_2}$.  

We claim that 
$K(A_1[\ell], A_2[\ell], Q/\ell) 
\subset K(A_1[\ell], A_2[\ell], P/\ell)$. 
If that was not the case, then $K(A_1[\ell], A_2[\ell], P/\ell, Q/\ell)$ would
have degree at least $\ell$ over   $K(A_1[\ell], A_2[\ell], P/\ell)$.
Therefore the density of
the set of primes in $K(A_1[\ell], A_2[\ell])$ which split
completely in both $K(A_1[\ell], A_2[\ell], P/\ell)$, and
$K(A_1[\ell], A_2[\ell], Q/\ell)$ will be at most  $1/\ell$ of the density of
the set of primes in $K(A_1[\ell], A_2[\ell])$ which split
completely in $K(A_1[\ell], A_2[\ell], P/\ell)$. 
However, we will prove below that the set of primes in
$K(A_1[\ell], A_2[\ell])$ which split in $K(A_1[\ell], A_2[\ell], P/\ell)$
and in $K(A_1[\ell], A_2[\ell], Q/\ell)$ is of much larger density, proving
that $K(A_1[\ell], A_2[\ell], Q/\ell)$ is a subfield of 
$ K(A_1[\ell], A_2[\ell], P/\ell)$.

A prime in $K(A_1[\ell], A_2[\ell])$ splits in $K(A_1[\ell], A_2[\ell], 
P/\ell)$
if $P$ has an $\ell$-division point  in the corresponding residue field
of $K(A_1[\ell], A_2[\ell])$. Now there are two possibilities under which
$P$ will have an $\ell$-division point. One, in which the order of $P$
in the residue field is coprime to $\ell$, in which case $P$ automatically
will have an $\ell$-division point. Second, in which $P$ has order
divisible by $\ell$ in which case if $P$ has an $\ell$-division point
in the residue field of $K(A_1[\ell], A_2[\ell])$, then clearly there
will exist an $\ell^2$-torsion point of $A_1$ which is also defined over
this residue field of $K(A_1[\ell], A_2[\ell])$. 

The field $K(A_1[\ell^2], A_2[\ell])$ is an abelian extension of
$K(A_1[\ell], A_2[\ell])$ with Galois group isomorphic to 
$Mp_{2g_1}({\Bbb Z}/\ell)$ which is the kernel of the natural map from 
$GSp_{2g_1}({\Bbb Z}/\ell^2)$ to $GSp_{2g_1}({\Bbb Z}/\ell)$.
For any point $Z$ belonging to $A_1[\ell]$,
the field $K(A_1[\ell], A_2[\ell], Z/\ell)$ is an abelian
extension of $K(A_1[\ell], A_2[\ell])$ of degree $\ell^{2g_1}$. 
Each of these extensions correspond (by Galois theory)
to subgroups $H_Z$  of $Mp_{2g_1}({\Bbb Z}/\ell)$.  
The group $Mp_{2g_1}({\Bbb Z}/\ell)$ operates 
on $({\Bbb Z}/\ell)^{2g_1}$, and $H_Z$ is the
stabilizer of the point $Z$ in $({\Bbb Z}/\ell)^{2g_1}$. This
implies in particular that every element of $H_Z$ has 1 as an
eigenvalue. Thus the determinant of $(1-X)$ is 0 for all $X$
belonging to $H_Z$, and for any choice of $Z$. Since the determinant
$\det(1-X)$ is a nonzero polynomial function on $Mp_{2g_1}({\Bbb Z}/\ell)$, 
the set of zeros is a proper subset of cardinality bounded by a 
(fixed) polynomial $P_{d-1}(\ell)$ 
in $\ell$ of degree $d-1$ where $d$ is the integer
such that $\ell^d$ is the cardinality of $Mp_{2g_1}({\Bbb Z}/\ell)$.

The set of primes of $K(A_1[\ell], A_2[\ell])$ for which
$A_1$ has a point of order $\ell^2$ in its residue field
correspond to primes $v$ of $K(A_1[\ell], A_2[\ell])$ such 
that $v$ splits in at least  one of the extensions 
$K(A_1[\ell],A_2[\ell], Z/{\ell}) \subset 
K(A_1[\ell^2], A_2[\ell])$. Thus the corresponding Artin symbol 
for the field extension $K(A_1[\ell^2], A_2[\ell])$ of
$K(A_1[\ell], A_2[\ell])$ must belong to one of the subgroups
$H_Z$ of $M_{2g_1}({\Bbb Z}/\ell)$. 

The union 
of the subgroups  
$H_Z$ has order bounded by $P_{d-1}(\ell)$. Thus the density 
of the set of primes in $K(A_1[\ell], A_2[\ell])$ which are
split in $K(A_1[\ell], A_2[\ell], P/\ell)$ but in none of
$K(A_1[\ell], A_2[\ell], Z/{\ell})$ is of density at least
equal to $[\ell^d-P_{d-1}(\ell)]/\ell^{2g_1+d}$ which is much
larger than $1/\ell^{2g_1+1}$ (for large $\ell$). For these primes, 
the order of $P$
in the residue field is coprime to $\ell$, and since order
of $Q$ divides the order of $P$, the order of $Q$ too
is coprime to $\ell$ in the residue field. Thus these primes
also split in $K(A_1[\ell], A_2[\ell], Q/\ell)$, and are of density greater 
than $1/{\ell^{2g_1+1}}$, proving our earlier claim that
$K(A_1[\ell], A_2[\ell], Q/\ell) \subset K(A_1[\ell], A_2[\ell], P/\ell)$. 
Hence by theorem 4 (due to [JR])
 applied to the abelian variety $A=A_1 \times A_2$, and points 
$(P,0),(0,Q)$, we get an endomorphism of $A$ taking 
$(P,0)$ to an integral multiple of $(0,Q)$. This implies that there is an
isogeny $j$ from $A_1$ to $A_2$ such that  $j(P)=nQ$ for some 
integer $n$, proving the theorem.

\section{ Proof of Theorem 2 (part 1)} 

To avoid trivialities, we will assume in what follows that the image of $A_1(K)$ under the homomprhism being considered is
infinite.

Given an abstract  homomorphism $$\phi:A_1(K) \rightarrow A_2(K)$$ that 
specialises mod $v$ for almost all places $v$ of $K$, it is clear that 
for any point $P$ of infinite order in $A_1(K)$, the order of $\phi(P) $ $\mod v$ which is a 
point in $A_2(k_v)$ 
divides the order of $P$ $\mod v$ which is a point in $A_1(k_v)$. Thus by theorem 1, there is an
isogeny $j$ and an integer $n$ (both depending on the point $P$) 
such that $j(P) = n\phi(P)$. What we need to prove for the part 1 of theorem 2
is that $j$ and $n$ can be chosen independent of $P$, say on a torsion-free
subgroup of $A_1(K)$. Note that if we can achieve $j(P_i)= n\phi(P_i)$ for the same $n$
(but perhaps different $j$)
for a set of points $P_i$, we can assume that that $n$ works for the subgroup
of $A_1(K)$ generated by $P_i$. Since 
the equation $j(P) = n\phi(P)$  can be multiplied by any integer, we can thus 
assume that
$n$ is independent of $P$. Thus for each torsion-free point $P$ of $A_1(K)$, we assume now 
that there is an isogeny $j$ dependent on $P$, and an integer $n \not = 0$ independent of $P$, 
 such that $j(P) = n \phi(P)$. 
The following general lemma implies that $j$ is independent of $P$. (The proof of this
fact, and of the lemma, for the case when End$(A_1) = {\Bbb Z}$ is rather trivial;
the lemma has been formulated in this generality because of the possible application
to the more general context.)

\begin{lemma}\label{abstract}
Let $A$ be a finitely generated free abelian group. Let $
{\cal D}$ be a division algebra which contains ${\Bbb Q}$ and is finite
dimensional over ${\Bbb Q}$. Let ${\cal O}$ be an order in ${\cal D}$. 
Suppose ${\cal O}$ acts on $A$ on the left 
making it into a left ${\cal O}$-module. Suppose that $f$ is an endomorphism
of $A$ as an additive group 
such that for all $a \in A$, there exists $f_a \in {\cal O}$ 
such that $f(a) = f_a \cdot a$. Then $f$ is multiplication by an element
of ${\cal O}$. 
\end{lemma}

\noindent{\bf Proof :} Clearly $A \otimes_{\Bbb Z} {\Bbb Q}$ is a vector
space over ${\Bbb Q}$ on which ${\cal O}\otimes_{\Bbb Z} {\Bbb Q} = {\cal D}$
acts, making it into a ${\cal D}$-vector space. From the hypothesis that
$f(a) = f_a \cdot a$ for all $a \in A$, we find that any ${\cal D}$-subspace of $A
\otimes_{\Bbb Z} {\Bbb Q}$ is stable under $f$ (extended to $A
\otimes_{\Bbb Z} {\Bbb Q}$). Write $A\otimes_{\Bbb Z} {\Bbb Q} = L_1 \oplus 
\cdots \oplus L_n$, as a direct sum of ${\cal D}$-subspaces 
$L_i$ of dimension
1 which as has been noted is invariant under $f$, i.e., $f(L_i) \subset
L_i$. Write $M_i =  L_i \cap A$; then $M_i$ is a lattice in $L_i$, and
$\oplus M_i$ is a subgroup of finite index in $A$. Since both $L_i$ and $A$
are invariant under $f$, so is $M_i$. Also, each $L_i$ and $A$ being invariant
under ${\cal O}$, so is $M_i$. We will now prove that the restriction 
of $f$ to $M_i$ is given by multiplication by an element $f_i$ in ${\cal O}$.

We can clearly assume that $M_i$ is a lattice in ${\cal D}$ which is
invariant under ${\cal O}$. Also, after scaling by an element of ${\cal D}^*$
on the right, we can assume that the lattice $M_i$ in ${\cal D}$ contains
1. 

Suppose that $f_i(1) = \alpha_i \in {\cal O}$. We can thus after replacing 
$f_i$ by $f_i-\alpha_i$, assume that $f_i(1) = 0$. We would like to prove
that $f_i$ is identically zero. Assuming the contrary, let $x$ be an
element in $M_i \cap {\cal O}$ such that $f_i(x) \not = 0$. After scaling
$x$, we can moreover assume that $f_i(x)$ belongs to ${\cal O}$. 

By hypothesis, for every element $m \in {\Bbb Z}$,  there exists an element
$\lambda_m \in {\cal O}$ such that $f_i(m+x) = \lambda_m\cdot (m+x).$ For
an element $z \in {\cal D}$, let ${\rm Norm}(z)$ denote the determinant of
the left multiplication by $z$ on ${\cal D}$. Since $f_i(x) = f_i(m+x) =
\lambda_m\cdot(x+m)$, it follows that 
${\rm Norm}(m+x)$ and ${\rm Norm}(f_i(x)) $  
are integers in ${\Bbb Q}$, and ${\rm Norm}(m+x)$ divides
 ${\rm Norm}(f_i(x)) $. Since ${\rm Norm}(m+x)$ is a polynomial in $m$
with coefficients in  ${\Bbb Z}$ of degree equal to the dimension
of ${\cal D}$ over ${\Bbb Q}$ of leading term 1 and constant term ${\rm
Norm}(x)$, the polynomial ${\rm Norm}(m+x)$ as $m$ varies takes arbitrary 
large values, hence cannot divide the fixed integer ${\rm Norm}(f_i(x))$.

We have thus proved that $f$ restricted to any 1 dimensional ${\cal D}$
submodule of 
$A\otimes_{\Bbb Z}{\Bbb Q}$ is multiplication by an  element of ${\cal D}$.    
From this it is trivial to see that the action of $f$ on 
$A\otimes_{\Bbb Z}{\Bbb Q}$ is multiplication by an  element of ${\cal D}$,
which must moreover lie in ${\cal O}$, completing the proof of the lemma.

\vspace{4mm}
\noindent{\bf Remark :} The lemma above holds good only in the integral
version stated above, and not for vector spaces, and hence is not totally
trivial. We point out an example to illustrate that the analogue
of the lemma is not true for vector spaces. For this, let $K$ be a finite
extension of ${\Bbb Q}$ of degree $> 1$. There is an action of $K^*$ on $K$ via left or
right multiplication. Let $f$ be an automorphism of $K$ (considered as a 
vector space over ${\Bbb Q}$) which does not
arise from the action of an element of $K^*$. Such an automorphism $f$ 
satisfies the hypothesis of the previous lemma as for any  $a \neq 0$, 
$f(a) \in K^*$, hence $f(a) = f_a \cdot a$ with $f_a \in K^*$, but such an
$f$ does not satisfy the conclusion of the lemma.

\section{Proof of theorem 2 (part 2)}
The proof of the 2nd part of theorem 2 for CM elliptic curves follows
from the theorem 2 of [CS] when combined with lemma 1. We will thus concentrate
our efforts in proving this part under the assumption that  End$_{\bar K}(A) 
= {\Bbb Z}$.

We begin with some preliminary results.

The following lemma is well known, cf. S.Lang's book, Algebra, section 10 of 
the chapter on Galois theory for $i=1$. It also follows from generalities on 
cohomology once we know it is true for $i=0$, which is of course clear.

\begin{lemma} Let $G$ be a group, and $E$ a $G$-module. Let $\tau$ be an 
element in the center of $G$. Then $H^i(G,E), i=0,1,\ldots$ is annihilated
by the map on $H^i(G,E)$ induced from the map $x \rightarrow \tau x-x$ from 
$E$ to itself.
\end{lemma}

\begin{lemma} 

Let $A$ be an abelian variety over a number field $K$. 
Let $K_{\ell^n} = K(A[\ell^n])$, and $G_{\ell^n}= {\rm Gal}(K_{\ell^n}/K)$. 
 Let $G_{\ell^\infty} = {\rm Gal}(K_{\ell^\infty}/K)$ 
with $K_{\ell^\infty} = \cup_{n}K_{\ell^n}$.
Then $H^1(G_{\ell^m},A[\ell^n])$ $(m \geq n$) 
and $H^1(G_{\ell^{\infty}},A[\ell^n]) $ 
are of finite
orders, bounded independent of $m$ and $n$. 
\end{lemma}

\noindent{\bf Proof :} We note that $G_{\ell^\infty}$ 
being a compact  $\ell$-adic Lie group,
is topologically finitely generated, hence each finite quotient such as 
$G_{\ell^m}$ is generated by a set of elements of cardinality independent of $m$. 

From the definition of $H^1(G,A)$ in terms of maps $\phi$ from
$G$ to $A$ such that $\phi(g_1g_2) = \phi(g_1)+g_1\phi(g_2)$, it follows that an element
of $H^1(G,A)$ is determined by a map on a set of generators of $G$. Since $A[\ell^n] \cong
({\Bbb Z}/\ell^n)^{2g}$ as abelian groups, it follows that $H^1(G_{\ell^m},A[\ell^n])$ 
is a finitely generated abelian group which is generated by a set of elements of
cardinality independent of $m,n$.    

It follows from a theorem of Bogomolov, cf. [B], that the $\ell$-adic Lie group 
$G_{\ell^\infty}$ contains homotheties congruent to 1 modulo $\ell^N$ for some integer
$N > 0$. Therefore by Lemma 1, 
$H^1(G_{\ell^m},A[\ell^n])$ is annihilated by $\ell^N$.  
It follows that $H^1(G_{\ell^m},A[\ell^n])$ 
is a finitely generated abelian group which is generated by a set of elements of
cardinality independent of $m,n$, and annihilated by $\ell^N$, and thus is
 of finite order,
bounded independent of $m,n$. 

The statement about $H^1(G_{\ell^{\infty}},A[\ell^n]) $ follows either by noting that
the cohomology $H^1(G_{\ell^{\infty}},A[\ell^n]) $ can be calculated
 in terms of continuous cochains 
on $G_{\ell^\infty}$, for which the earlier argument applies as well, or by noting that
$H^1(G_{\ell^{\infty}},A[\ell^n]) $ is the direct limit of $H^1(
G_{\ell^{m}},A[\ell^n]) $ (direct limit over $m$), and a direct limit of  
finitely generated abelian groups  each of which is generated by a set of elements of
cardinality independent of $n$, and each annihilated by $\ell^N$, is of order 
bounded independent of $n$. 

\begin{lemma}\label{div}
  Given an abelian variety $A$ over $K$, a point $P$  of $ A(K)$ of
  infinite order, and any prime
  $\ell$, there are infinitely many places $v$ of $K$  
(in fact of positive density) such that the reduction
  of $P$ mod $v$ has order divisible by $\ell$.
\end{lemma}

\noindent{\bf Proof:} We claim that there is a sufficiently large $n$
such that the extension
$K_{P,\ell^n}=K(A[\ell^n],{\frac{1} {\ell^n}}.P)$ is a non-trivial
extension of  $K_{\ell^n}=K(A[\ell^n])$. If we grant the claim
then there is a positive density of $v$ that split in $K(A[\ell^n])$ (we
denote
the Galois group ${\rm Gal}(K_{\ell^n}/K)$ by $G_{\ell^n}$) but not in
$K(A[\ell^n],{\frac{1} {\ell^n}}.P)$, and by inspection we see that for
such $v$'s
the reduction of $P$ mod $v$ has order divisible by $\ell$.

From Lemma 2, it follows that $H^1(G_{\ell^n},A[\ell^n])$
is bounded independently of $n$. Using maps between the Kummer sequences
\begin{equation*}
\begin{CD}
0 @> >>  A(K)/\ell^nA(K) @> >> H^1(G_K,A[\ell^n])
@> >> H^1(G_K,A)[\ell^n]  @> >> 0
\\ 
&& @V VV @V VV @V VV 
\\
0 @> >>   A(K_{\ell^n})/\ell^nA(K_{\ell^n}) @> >>
H^1(G_{K_{\ell^n}},A[\ell^n]) @> >>    H^1(G_{K_{\ell^n}},A)[\ell^n]
@>  >>     0
\end{CD}
\end{equation*}
the claim follows for large enough $n$, putting together the observation that
the kernel of the restriction map 
$H^1(G_K,A[\ell^n]) \rightarrow H^1(G_{K_{\ell^n}},A[\ell^n])$ 
is $H^1(G_{\ell^n},A[\ell^n])$, and the
fact that as $P$ is non-torsion
the image of $P$ under the
coboundary map (in the first exact sequence) in $H^1(G_K,A[\ell^n])$ has unbounded order as $n$
varies. From this the claim follows and hence the lemma.

\vspace{3mm}

We next have the following lemma.

\begin{lemma}\label{notdiv}
  Given an abelian variety  $A$ over a number field $K$  satisfying the
 hypotheses of 
Theorem 2, a point $P$ of $ A(K)$ of
  infinite order, and any prime
  $\ell$, there are infinitely many places $v$ of $K$
 (in fact of positive density) such that the reduction
  of $P$ mod $v$ has order prime to $\ell$.
\end{lemma}

\noindent{\bf Proof:} 
From the proof of Lemma 3, it follows that $K_{P,\ell^{\infty}}$
is an infinite extension of $K_{\ell^\infty}$. Under the hypothesis of this
lemma, $G_{\ell^\infty}$ is an open subgroup of finite index of $GSp_{2g}(
{\Bbb Z}_\ell)$ by a theorem of Serre [S-IV]. (If the abelian variety is not
principally polarised, the image of the Galois group does not necessarily 
 sit inside a symplectic similitude group, and hence a slight 
modification to the argument given below 
will have to be made which we leave to the reader.)
 If we denote by
$E_{\ell^\infty}
$ the
Galois group of $K_{P,\ell^{\infty}}$ over $K$, and $A_{\ell^{\infty}}$,
the Galois group of $K_{P,\ell^{\infty}}$ over $K_{\ell^{\infty}}$, we have
an exact sequence of groups,
$$ 0 \rightarrow A_{\ell^\infty}
\rightarrow E_{\ell^\infty}
 \rightarrow G_{\ell^\infty} \rightarrow 1.$$
This exact sequence is a subsequence 
 of the following
(split) exact sequence of ${\Bbb Z}_\ell$-rational points of algebraic
groups:
$$ 0 \rightarrow {\Bbb Z}_{\ell}^{2g} \rightarrow E({\Bbb Z}_\ell)
 \rightarrow GSp(2g,{\Bbb Z}_{\ell})
\rightarrow 1,$$
(i.e., 
$A_{\ell^\infty} \subset {\Bbb Z}_{\ell}^{2g}, 
 E_{\ell^\infty} \subset E({\Bbb Z}_\ell),
 G_{\ell^\infty} \subset GSp(2g,{\Bbb Z}_{\ell})$),
in which $E({\Bbb Z}_\ell)$ is the semi-direct product 
$GSp(2g,{\Bbb Z}_{\ell}) $ 
with ${\Bbb Z}_{\ell}^{2g}$. The embedding of $E_{\ell^{\infty}}$ inside 
$E({\Bbb Z}_\ell)$ is obtained by choosing a sequence of points $P_n$ 
in $A$ with $\ell \cdot P_1 =P$, and $\ell \cdot P_{i+1} = P_i$. 
Since the action of $GSp(2g,{\Bbb Z}_{\ell})$ on  ${\Bbb Z}_{\ell}^{2g}$ is
irreducible, and $G_{\ell^\infty}$ is an open subgroup of
$GSp(2g,{\Bbb Z}_{\ell}) $, it follows that $A_{\ell^\infty} $ is an
open subgroup of ${\Bbb Z}_{\ell}^{2g}$, and hence $E_{\ell^\infty}$
is an open subgroup of the
semi-direct product $GSp(2g,{\Bbb Z}_{\ell}) $ with
${\Bbb Z}_{\ell}^{2g}$.

We will prove that the intersection of the fields $K_{P,\ell^n}$ and
$K_{\ell^{n+1}}$ is $K_{\ell^n}$ for some $n$ large enough.
This, together with the theorem of
Bogomolov recalled earlier, will imply that there is a positive
density of primes in $K$ which are split in $K_{P,\ell^n}$ and for which
the Frobenius as an element of ${\rm Gal }(K_{\ell^{n+1}}/K)$ is a
non-trivial homothety in $GSp(2g,{\Bbb Z}/\ell^{n+1})$ which is congruent to
1 modulo $\ell^n$. For such primes $v$, it can be easily seen that
the order of $P$ modulo $v$ is not divisible by $\ell$,
completing the proof of the lemma.
 
It thus suffices to prove that the intersection of the fields
$K_{P,\ell^n}$ and  $K_{\ell^{n+1}}$ is $K_{\ell^n}$.
 
Let $E_{\ell^n}$ (resp. $G_{\ell^n}$) denote the Galois group of
$K_{P,\ell^n}$
(resp. $K_{\ell^n}$) over $K$, and let $A_{\ell^n}$ denote the Galois group
of $K_{P,\ell^n}$ over $K_{\ell^n}$. We have the exact sequence of groups,
$$ 0 \rightarrow A_{\ell^n} \rightarrow E_{\ell^n} \rightarrow G_{\ell^n}
 \rightarrow 1.$$
 
It is clear  that the intersection of the fields
$K_{P,\ell^n}$ and  $K_{\ell^{n+1}}$ is $K_{\ell^n}$ if and only if
inside the group $E_{\ell^n}$ which is a quotient of $E_{\ell^{n+1}}$,
the image
of $A_{\ell^{n+1}}$ is $A_{\ell^n}$.
 
Since   $E_{\ell^\infty}$
 is an open subgroup of the
semi-direct product $GSp(2g,{\Bbb Z}_{\ell}) $ with
${\Bbb Z}_{\ell}^{2g}$, it follows that $E_{\ell^\infty}$ contains
the natural congruence subgroup of level $\ell^n$ in this semi-direct
product
for some $n$, say $n=n_o$. From this it is clear that $E_{\ell^{n+1}}$ is
the
full inverse image of $E_{\ell^n}$ under the natural mapping from
$E({\Bbb Z}/\ell^{n+1})$ to $E({\Bbb Z}/\ell^{n})$, $n \geq n_0$,
 from which the surjectivity
of the mapping from $A_{\ell^{n+1}}$ onto $A_{\ell^n}$  clearly follows,
completing the proof of the lemma.

\vspace{3mm}

\noindent{\bf Remark:} 
The above proof can be generalised to yield
that for any abelian variety $A$ defined over $K$ and any point $P$ of 
$A(K)$  which does not project to a non-zero torsion point in 
any (geometric) subquotient of 
$A$, given a prime $\ell$ there are a positive density of 
places $v$ of $K$ such that $P$ mod $v$ has order prime to $\ell$.

\vspace{3mm}

\noindent{\bf Remark:} Although we have given separate proofs of 
Lemmas 4 and 5, observe that Lemma 5 implies Lemma 4. This follows
by applying Lemma 5 to the point of infinite order $P+R$ where $R$ is a 
nonzero $\ell$-torsion point. (We note that to prove lemma 4, we are allowed
to go to a finite extension of $K$, and hence assume that $A$ has 
nonzero $\ell$-torsion point over $K$.) 
Lemma 5 implies that there are infinitely many places 
$v$ of $K$ for which the order of $P+R$ is coprime to $\ell$. It is easy to
see that for such places, the order of $P$ must be divisible by $\ell$.

\begin{cor}\label{triv}
Under the hypothesis and notation of theorem 2, the image of any 
non-torsion point $P$ in $A_1$ is non-torsion in $A_2$. 
\end{cor}

\noindent{\bf Proof:} Immediate from Lemma \ref{notdiv}.

\vspace{3mm}

\noindent Lemma \ref{notdiv} allows us to strengthen theorem 1.

\begin{prop} Let $A_1$ and $A_2$ be two simple 
abelian varieties over a number field $K$ 
with the endomorphism rings over the algebraic closure of $K$ 
as ${\Bbb  Z}$.  Assume that  ${\rm dim}(A_i)=g_i$ is 
either odd or $g_i=2$ or 6. (We do not assume that $\dim A_1 = \dim A_2$.)
Let $P$ (resp. $Q$) be a  non-torsion point in $ A_1(K)$ (resp. $A_2(K)$) 
such that for almost
  all places $v$ of $K$, the order of $Q$ mod $v$ divides the order of $P$
  mod $v$. Then $A_1$ and $A_2$ are isogenous, and 
there is  an isogeny $j_0$ from $A_1$ to $A_2$ such that $j_0(P) = Q$. 
Furthermore, if one of the abelian
varieties is an elliptic curve, we have the same conclusion without
any restriction on its endomorphism ring.
\end{prop}

\noindent{\bf Proof:} From theorem 1, 
there is  an isogeny $j$ from $A_1$ to $A_2$ such that $j(P) = nQ$ for some
integer $n$. We wish to prove that $n$ can be chosen to be 1. 
Let $\ell^r$ be the highest power of a prime $\ell$ which divides 
$n$. We will prove that $\ell^r$ torsion points of $A_1$ are contained
in the kernel of $j$. 
If that were not the case, there would be an $\ell^r$
torsion point of $A_1$, say $R$, which does not belong to the kernel of $j$.
By lemma 5, there are infinitely many places $v$ of $K$ such that the order 
of $P+R$ is coprime to $\ell$ in the residue field $k_v$ of $K$, and hence the
$\ell$-primary component of the orders of $P$ and $R$ are the same, 
of order $\ell^r$. For such
places $v$, the order of $j(P+R)=j(P)+j(R)= nQ+j(R) $ is also 
coprime to $\ell$. By choice, $j(R)$ is a nonzero torsion point on
$A_2$ of order dividing $\ell^r$, say $\ell^s, 0< s \leq r$. Thus since the order of $nQ+j(R)$ 
is coprime to $\ell$, the $\ell$-primary components of the order of $j(R)$ 
and $nQ$ are the same, hence the $\ell$-primary component of the order of $nQ$ is
$\ell^s$, and therefore of $Q$, $\ell^{r+s}$, contradicting our hypothesis that
order of $Q$ divides the order $P$ at each place. This 
proves that $\ell^r$ torsion points of $A_1$ are contained in the kernel of $j$ where
$\ell^r$ is the highest power of $\ell$ dividing $n$. Thus all the $n$-torsion 
points of $A_1$ is contained in the kernel of $j$.  Therefore the isogeny $j$ 
can be written as $nj_0$ for an isogeny $j_0$ from $A_1$ to $A_2$.
Thus we have $n(j_0(P)-Q)=0$. This implies that $j_0(P) = Q +S$ for a certain torsion
point $S$ on $A_2$. 

If $S$ is nonzero, let the order of $S$ be divisible
by a prime $\ell$.  Choose a place $v$ of $K$ where the order of $P$, and hence
of $Q$ and $j_0(P)$ are coprime to $\ell$. However, since $j_0(P) = Q +S$, its 
order is divisible by $\ell$, a contradiction to $S$ being nonzero.

\vspace{4mm}

\noindent {\bf End of the proof of part 2 of Theorem 2.}
Choose a torsion-free subgroup $B$ of $A_1(K)$ such that
$A_1(K)=A_1(K)_{tors} \oplus B$. Then  
for each $P \in B$
we see from Corollary \ref{triv} that $\phi(P)$ is not torsion, and
by Proposition 1 that  $\phi(P)=j_0(P)$ for some isogeny $j_0$ that might a priori
depend on $P$.  Considering $P+Q$ for $P,Q \in B$ it follows that $j_0$ is independent
of $P$. Now considering $\phi(P+P')=j'_0(P+P')$ and $\phi(P)=j_0(P)$ where 
$P'$ is in $ A(K)_{tors}$ and $P \in B$, we conclude that $\phi(P')-j'_0(P') =
(j_0-j'_0)(P) 
$. The left hand side of this equation is a point of finite order, whereas 
right hand is of infinite order unless $j_0=j'_0$, forcing $j_0 = j'_0$, and 
thus we are done with the proof of Theorem 2.

\section{Rigidity for arithmetic groups}

We begin with the following theorem for tori.

\begin{prop}
  Given homomorphism $\phi:{\cal O}_K^* \rightarrow {\cal O}_K^*$ that
  reduces mod $v$ for almost all places $v$ of a number field $K$,
  then $\phi$ is induced by the $m$th power map for some integer $m$.
\end{prop}

\noindent{\bf Proof:} The proof is a direct consequence of Theorem 1
of [RS] and the fact that any finite subgroup of $K^*$ is cyclic.

\vspace{3mm}

We next have the following theorem for arithmetic groups.

\begin{theorem} Let $\Gamma$ be a subgroup of $SL(2,{\Bbb Z})$ of finite 
index. Let $\phi$ be a non-trivial homomorphism of $\Gamma$ into itself. Assume
that for all primes $p$ in an infinite set $S$ of primes, 
$\phi$ factors to give a homomorphism
$\phi_p: SL(2,{\Bbb Z}/p) \rightarrow SL(2,{\Bbb Z}/p)$
$$
\begin{CD}
  \Gamma @>\phi>> \Gamma \\
  @V    VV    @VV    V \\
  SL(2,{\Bbb Z}/p) @>\phi_p>>  SL(2,{\Bbb   Z}/p).
\end{CD}
$$
 Then $\phi$ is an automorphism which is the restriction to $\Gamma$ of the
inner-conjugation action of an element in $GL(2,{\Bbb Q})$.
\end{theorem}

\noindent{\bf Proof:} Let $ A$ be the ring which is the direct product 
of ${\Bbb Z}/p$ for all $p$  in $S$.
Clearly ${\Bbb Z}$ is a subring of $A$, and there is thus
an injective  homomorphism from $SL(2,{\Bbb Z})$ to $SL(2,A)$. 
Since there is an injective  homomorphism from $SL(2,{\Bbb Z})$ to 
$SL(2,A)$ for $A$, the direct product of {\it any} infinite set of primes, 
it is clear that  $\phi_p$ can be trivial for at most 
finitely many $p$ in $S$.  After 
replacing $S$ by this slightly smaller set, we assume that $\phi_p$ is 
surjective for all $p$ in $S$, and hence the $\phi_p$ are given by the 
inner-conjugation action of an element $g_p $ in $GL(2,{\Bbb Z}/p)$. Here
we are using the well-known facts:

\begin{enumerate}
\item  any surjective 
homomorphism of $SL(2,{\Bbb Z}/p)$
into itself is given by the inner-conjugation of an element of 
$GL(2,{\Bbb Z}/p)$. 
\item  any  
homomorphism of $SL(2,{\Bbb Z}/p)$
into itself is either trivial or is surjective if $p>3$.
\end{enumerate}
 
From this we see that the representations $\phi: \Gamma \rightarrow GL_2({\Bbb Q})$
and the ``identity'' representation ${ id}:\Gamma  \rightarrow GL_2({\Bbb Q})$ have the
same trace. Further the second representation is irreducible.
From this we conclude that $\phi$ and ${ id}$ are conjugate
by an element of $GL_2({\Bbb Q})$ and further as $\phi$ goes into $\Gamma$
by comparing covolumes we see that $\phi$ is an automorphism of $\Gamma$.

\vspace{4mm}

\noindent{\bf Remarks: } 1. The above proof is due to Serre: we had a different proof 
in an earlier version.

2. To deduce rigidity results for $SL_2$ one just 
needs that the abstract homomorphism specialises for any infinite
set of primes rather than for almost all or even a positive density of
primes which is crucial for abelian varieties. 
Further unlike
abelian varieties the analog of Question 1 for $SL_2$ dealt
with here is not a ``pointwise question'' (see remark following the question).

3. The proof works for $SL(n,{\Bbb Z})$ for  any $n$ to say that if $\phi$
is a homomorphism of a subgroup of finite index of $SL(n,{\Bbb Z})$ 
onto another subgroup of finite index 
$SL(n,{\Bbb Z})$   which specialises for infinitely
many primes $p$ to give a homomorphism of $SL(n,{\Bbb Z}/p)$ to itself 
(we recall that an automorphism of
$SL(n,{\Bbb Z}/p)$   is generated by inner automorphism from
$GL(n,{\Bbb Z}/p)$, and the automorphism $A \rightarrow {}^tA^{-1}$), then
$\phi$ is algebraic.

4. Because of the strong rigidity theorem, the algebraicity of abstract
homomorphisms of arithmetic lattices in semi-simple Lie groups is of
interest only for arithmetic lattices in $SL(2,{\Bbb R})$ which are constructed
using division algebras over totally real number fields. Our method clearly
applies for such lattices too.

5. The following theorem when combined with the previous one, 
proves that for a subgroup $\Gamma$ of finite index 
in $SL(2,{\Bbb Z})$, any homomorphism of $\Gamma$ into itself which extends
to $SL(2,A)$ for $A = \prod_{p \in T} {\Bbb Z}/p$, $T$ an infinite set 
of primes, must be given by an inner-conjugation action by an element of $GL(2,
{\Bbb Q})$, giving a different perspective to the earlier theorem.

\begin{theorem} Let $G$ be any simply-connected, split, semisimple 
algebraic group
over ${\Bbb Q}$.  Then any homomorphism of $G(A)$ to itself where
$A$ is the
 product of ${\Bbb Z}/p$ for primes $p$ belonging to a set $T$ that may be 
finite or infinite
 (and contains only sufficiently large primes) 
is factorisable, i.e., any homomorphism
$f$ from $G(A)$ to itself is of the form $\prod_{p \in A} f_p$ for certain
homomorphisms $f_p$ from $G({\Bbb Z}/p)$ to itself.
 
\end{theorem} 
 
\noindent{\bf Proof: }  We will accomplish the proof of this theorem in 
several steps.
 
\begin{enumerate}
\item Any homomorphism from $G({\Bbb Z}/p)$ to $G({\Bbb Z}/q)$, 
$p$ not $q$ is trivial (for $p$ a sufficiently large prime).

 Assume that the mapping is non-trivial. Then as $G({\Bbb Z}/p)$ is a simple 
group modulo its center,
any homomorphism from $G({\Bbb Z}/p)$ to $G({\Bbb Z}/q)$
must be injective when restricted to unipotent elements in $G({\Bbb Z}/p)$.
 Because $p$ is not $q$, image of a unipotent element in 
$G({\Bbb Z}/p)$ cannot have a unipotent component
 in the Jordan decomposition in $G({\Bbb Z}/q)$.
 So image of any unipotent element in $G({\Bbb Z}/p)$ is semi-simple in
$G({\Bbb Z}/q)$. 

Note that a unipotent in $SL_2({\Bbb Z}/p)$
 has many powers that are conjugate to itself. By Jacobson-Morozov 
(which is applicable since we are looking only at large primes), 
the same holds good about non-trivial unipotents in $G({\Bbb Z}/p)$. 
Hence the  image of a non-trivial unipotent in $G({\Bbb Z}/p)$ too will have
 many distinct powers that are conjugate to itself. But a semi-simple
 element in $G({\Bbb Z}/q)$  has at most $|W|$ many
powers that are conjugate to itself, where $|W|$ denotes the order of the
Weyl group of $G$, completing the proof of this step.

\item Step 1 proves that any homomorphism from  $G(A)$ to itself 
 when restricted to direct sum is factorisable. However going from
 direct sum to direct product needs some more arguments and essentially
 the following step suffices.
 
\item Any homomorphism  from  $G(A^S)$ to  $G({\Bbb Z}/p)$ 
 must be trivial for some finite set $S$ (depending on $p$) of primes in 
 $A$ with
 $A^S=A-S$ where $S$ is the set of prime divisors in $A$ of the cardinality of
 $G({\Bbb Z}/p)$, which we denote by $d$.
 
 We first prove that a unipotent in $G(A^S)$ must go trivially to
 $G({\Bbb Z}/p)$. Again by Jacobson-Morozov (applied to the ring
$A^S$ which is a product of fields), 
this would follow if we can prove that under any homomorphism from $SL_2(A^S)$ 
to $G({\Bbb Z}/p)$, any unipotent in $SL_2(A^S)$ must go trivially in 
$G({\Bbb Z}/p)$.
But this follows because multiplication by $d$ is an
 isomorphism on $A^S$ whereas it is trivial on $G({\Bbb Z}/p)$.
 
 We will be done if we can prove that unipotents in $G(A)$ for any set of
 primes $A$ generates  $G(A)$. (This step is not true for an arbitrary
ring $A$, but is true here as $A$ is a product of fields.)

\item Let $B$ be any Borel subgroup in $G$ defined over ${\Bbb Q}$. We 
will prove that any element of $B(A)$ belongs to the group generated by
the  unipotents. 

We need to prove that for a torus $T$ contained in $B$, and defined over 
${\Bbb Q}$, elements of $T(A)$ belong to the subgroup generated by the 
unipotents in $G(A)$. 
 
For $SL_2(A)$, this follows from the matrix identity 
which is a product of 4 unipotent
 matrices taken from Deligne's article in Modular Forms vol 2, Springer 
Lecture Notes in Mathematics, vol. 349.

$$\left(\begin{matrix}
 a^{-1} & 0 \\ 
 0    &   a \end{matrix}\right ) 
= \left(\begin{matrix} 1 &  -a^{-1} \\    
                 0 &   1   \end{matrix}  \right )   
  \left (\begin{matrix} 1  &    0 \\
            a-1 &  1\end{matrix} \right )     
 \left (\begin{matrix} 1  &    1 \\
         0   &  1 \end{matrix}\right )     
 \left (\begin{matrix} 1  &    0 \\
            -(a-1)/a &  1\end{matrix} \right) $$

As $G$ is simply-connected, any element of $T(A)$ is a product of 
elements of 1 dimensional tori in $T(A)$ arising out of the image of
the diagonal torus in $SL_2$ under mappings of $SL_2$ to $G$ corresponding
to the simple roots.

\item  Unipotents in  $G(A)$ for any set of primes $A$ generates  $G(A)$.
    
This follows from step 4 combined with Bruhat decomposition 
(applied component-wise to write $g=(g_p)$ as $(u_pw_pb_p)$) and the fact that
any element of the Weyl group is generated by the unipotent elements.

\item It follows from step 1 and 3 that any homomorphism from $G(A)$ to
$ G({\Bbb Z}/p)$ is trivial on $G(A-p)$, hence we have proved that any 
homomorphism from $G(A)$ to itself is factorisable.

\end{enumerate}

\section{Rational Varieties}
After having considered the case of Abelian varieties and semi-simple 
algebraic groups, 
it is tempting to consider arbitrary maps of algebraic varieties 
over number fields which reduce
nicely under reduction modulo all finite primes, 
i.e. such that the following diagram
$$
\begin{CD}
  V(K) @>\phi>> V(K) \\
  @V{sp_v}    VV    @VV{sp_v}   V \\
  V(k_v) @>\phi_v>>  V(k_v)
\end{CD}
$$ 
commutes where $k_v$ is a residue field of the ring of integers of $K$,
 and to ask whether such maps come from an algebraic one on $V$. 
The question will
naturally be more meaningful if $V$ is assured of many rational points, 
if for 
instance, $V$ is a smooth projective rational variety.
It seems specially interesting to investigate it for
flag variety $G/P$ for a parabolic $P$ in a semi-simple split group 
$G$ over $K$. Here, we merely point out that the analogous
question for the affine line over ${\Bbb Z}$ is false, 
i.e. there exists a set-theoretic
map from ${\Bbb Z}$ to ${\Bbb Z}$ which is not polynomial but which makes
 the following diagram commute.
$$
\begin{CD}
  A^1({\Bbb Z}) @>\phi>> A^1({\Bbb Z}) \\
  @V{sp_p}    VV    @VV{sp_p}   V \\
  A^1({\Bbb Z}/p) @>\phi_p>>  A^1({\Bbb Z}/p).
\end{CD}
$$ 
This is constructed for instance using 
$$ \psi(n) = a_0+a_1n+a_2 {n(n-1)} +a_3{n(n-1)(n-2)} + \cdots,$$
an infinite sum, which reduces to a finite sum for each $n$, where
 $a_i$ are integral, and $n \geq 0$,  and defining $\phi(n) = \psi(n^2)$. We find  
that since $\phi(n)$ is congruent to $\phi(m)$ modulo any integer $N$ for which $m$ is
congruent to $n$ modulo $N$, the diagram above commutes for $\phi$.

\vspace{3mm}

\noindent{\bf Acknowledgement:} We thank J.-P.Serre for helpful correspondence,
and K. Ribet for pointing out the preprint [BGK] which has some overlap with the
present work. We note, however, that the theorems of [BGK] are 
all proved under the assumption that $A_1=A_2$, and that 
our finer theorem 2.2 is not considered there.

\section{References}

\noindent [BGK] Banaszak, G., Gajda, W., Kraso\'n, P., {\it A support
problem for the intermediate Jacobians of $\ell$-adic 
representations}, preprint.

\vspace{3mm}

\noindent [CS] Corrales, C., and Schoof, R., {\it The support
problem and its elliptic analogue}, J. of Number Theory 64 (1997), 276--290.

\vspace{3mm}

\noindent[B] Bogomolov, F., 
{\it Sur l'alg\'ebricit\'e des repr\'esentations $l$-adiques}, 
 C. R. Acad. Sci. Paris S\'er. A-B 290 (1980), no. 15, 701--703. 

\vspace{3mm}

\vspace{3mm}

\noindent [JR] Jacquinot, O., Ribet, K., {\it Deficient points
on extensions of abelian varieties by ${\Bbb G}_m$}, Journal of Number
Theory 25 (1987), 133--151.



\vspace{3mm}














\vspace{3mm}

\noindent [S]  Serre, J-P., {\it  Propri\'{e}t\'{e}s galoisiennes des
points d'ordre finides courbes elliptiques}, Invent. Math. 15 (1972), 259--331.

\vspace{3mm}

\noindent [S-IV] Serre, J-P., {\it Oeuvres}, 
Vol. IV, articles 133--138, Springer-Verlag, 2000.



\vspace{4mm}

\noindent{\it Addresses of the authors:} 

\vspace{2mm}

\noindent CK: Dept. of Math, University of Utah, 155 S 1400 E, Salt Lake City, UT 84112, USA: shekhar@math.utah.edu

\noindent School of Mathematics, Tata Institute of Fundamental Research, Colaba, Bombay-400005, INDIA. e-mail:
shekhar@math.tifr.res.in 
\vspace{2mm}

\noindent DP: Harish-Chandra Research Institute, Chhatnag Road, Jhusi,
Allahabad-211019, INDIA. e-mail: dprasad@mri.ernet.in

\end{document}